\documentclass[11pt,epsfig]{article}
\usepackage{mathrsfs}
\usepackage{amsmath}
\usepackage{amsthm}
\usepackage{epsfig}
\usepackage{graphicx}
\usepackage{hyperref}
\usepackage{amsfonts}
\usepackage{color}

\newtheoremstyle{thmm}{1.5ex plus 1ex minus .2ex}{1.5ex plus 1ex minus
.2ex}{\rmfamily}{}{\bfseries}{}{1em}{} \theoremstyle{thmm}
\newtheorem{theorem}{Theorem}[section]
\newtheorem{lemma}{Lemma}[section]

\renewenvironment{proof}[1][Proof]{\noindent\textit{#1. } }{\hfill$\square$}

\allowdisplaybreaks \textwidth  6.0in \textheight 9.0in \topmargin
-0.6in \oddsidemargin 0.2in \evensidemargin 0.0in

\newcommand{\nn}{\nonumber}
\def \endproof{\vrule height8pt width 5pt depth 0pt}
\def\refe#1{(\ref{#1})}

\def\R{\mathbb{R}}

\begin{document}

\title{\bf Unconditionally optimal error estimates
of a Crank--Nicolson Galerkin method for the\\
nonlinear thermistor equations
}

\author{Buyang Li\,\setcounter{footnote}{0}\footnote{Department
of Mathematics, City University of Hong Kong, Kowloon, Hong Kong.
The work of the
authors was supported in part by a grant from the Research Grants
Council of the Hong Kong Special Administrative Region, China
(Project No. CityU 102005) ~
{\tt libuyang@gmail.com} (B. Li), {\tt
hdgao2@student.cityu.edu.hk} (H. Gao), {\tt
maweiw@math.cityu.edu.hk} (W. Sun). },
~Huadong Gao\footnotemark[1]
~~and ~Weiwei Sun\footnotemark[1]}

\maketitle

\begin{abstract}
This paper focuses on unconditionally optimal error analysis
of an uncoupled and linearized Crank--Nicolson Galerkin finite
element method for the time-dependent nonlinear thermistor equations in
$d$-dimensional space, $d=2,3$.
We split the error function into two parts,
one from the spatial discretization and one from the temporal discretization,
by introducing a corresponding time-discrete (elliptic) system.
We present a rigorous analysis for the regularity of the solution of
the time-discrete system and error estimates of the time discretization.
With these estimates and the proved regularity, optimal error
estimates of the fully discrete Crank--Nicolson Galerkin method
are obtained unconditionally.
Numerical results confirm our analysis and show the efficiency of
the method.

\end{abstract}


\section{Introduction}
\setcounter{equation}{0}
We consider the time-dependent nonlinear thermistor system
\begin{align}
&\frac{\partial u}{\partial t}-\Delta u=\sigma(u)|\nabla\phi|^2,
\label{e-heat-1}
\\[3pt]
&-\nabla\cdot(\sigma(u)\nabla\phi)=0, \label{e-heat-2}
\end{align}
for $x\in\Omega$ and $t\in[0,T]$, where $\Omega$ is a
bounded domain in $\R^d$, $d=2,3$.
The initial and boundary conditions are
given by
\begin{align}
\label{BC}
\begin{array}{ll}
u(x,t)=0,\quad \phi(x,t)=g(x,t)~~
&\mbox{for}~~x\in\partial\Omega,~~t\in[0,T],\\[3pt]
u(x,0)=u_0(x)~~ &\mbox{for}~~x\in\Omega .
\end{array}
\end{align}

The nonlinear system above describes the model of electric heating
of a conducting body, where $u$ is the temperature, $\phi$ is the
electric potential, and $\sigma$ is the temperature-dependent
electric conductivity. Following the previous works \cite{EL,Zhao},
we assume that $\sigma\in W^{1,\infty}(\R)$ and
\begin{align}
\label{sigma} \sigma_1\leq \sigma(s)\leq \sigma_2,
\end{align}
for some positive constants $\sigma_1$ and $\sigma_2$.

Theoretical analysis for the time-dependent thermistor equations
was done by several authors \cite{AX,ALM,Cim,
Yuan2, YL}. Among these works, Yuan and Liu \cite{YL} proved the
existence and uniqueness of a $C^\alpha$ solution in three-dimensional space.
Based on their result, further regularity can be derived with suitable
assumptions on the initial and boundary conditions. Numerical methods
and analysis for the thermistor system can be found in \cite{AL,
AY, EL, Yue, Zhao,  ZW}. For the system in two-dimensional
space, the optimal $L^2$ error estimate of a mixed finite element method
with a linearized semi-implicit Euler scheme was obtained in
\cite{Zhao} under a weak time-step condition. Error analysis for the
three-dimensional model was given in \cite{EL}, in which a
linearized semi-implicit Euler scheme with a linear Galerkin FEM was
used. An optimal $L^2$-error estimate was obtained under the
condition $\tau = O(h^{1/2})$. A more general time discretization
with higher-order finite element approximations was studied in
\cite{AL}. An optimal $L^2$-norm error estimate was given under the
conditions $\tau = O(h^{3/2p})$ and $r \geq 2$, where $p$ is the
order of the time discretization and $r$ is the degree of piecewise
polynomials of the finite element space.

Clearly, there are several different time discretizations for
nonlinear parabolic systems, explicit, semi-explicit (or semi-implicit)
and implicit. The most popular and widely-used
approach is linearized (semi)-implicit scheme. At each time step,
the scheme only requires the solution of a linear system.
However, time-step condition is always a key issue for such a scheme.
To study the error
estimate of linearized (semi)-implicit schemes, the boundedness of
the numerical solution (or error function) in $L^{\infty}$ norm or a
stronger norm is often required. If a priori estimate for numerical
solution in such a norm cannot be provided, one may employ the
mathematical induction with an inverse inequality to bound the numerical
solution, such as, by
\begin{equation}
\| U_h^n - R_hu^n \|_{L^\infty} \le C h^{-d/2} \| U_h^n - R_hu^n \|_{L^2} \le
C h^{-d/2} (\tau^{m} + h^{r+1}) ,
\end{equation}
where $U^n_h$ is the finite element solution, $u^n$
is the exact solution and $R_h$ is certain projection operator.
The above approach, however, requires a time-step condition $\tau = O(h^{d/(2m)})$.
This approach has been widely used in the error analysis of many
different nonlinear parabolic PDEs, $e.g.$, see \cite{
He,KL,Liu1} for Navier-Stokes equations, \cite{AL, EL,
Zhao} for nonlinear thermistor problems, \cite{EW, SS,Wang} for
porous media flows, \cite{CL,WHS} for viscoelastic fluid
flow, \cite{MS} for KdV equations, \cite{chen,MH} for
the Ginzburg-Landau equations, \cite{BC,Tou} for
nonlinear Schr\"{o}dinger equations and \cite{DM,WML} for some
other equations. In all these works, error estimates were
established under certain time step restrictions.
The time-step restrictions arising from theoretical analysis
may result in the use of a very small time step
and extremely time-consuming in practical computations.
However, we believe that such time-step conditions may not be necessary
for most cases. A new approach was introduced in
our recent works \cite{LS1, LS2}, also see \cite{Li},
in which the error estimates of a linearized
backward Euler Galerkin methods for a porous media flow and
the thermistor system were obtained, respectively, under the condition
of $h$ and $\tau$ being smaller than a positive constant.
In this paper, we propose an uncoupled and linearized Crank--Nicolson Galerkin
finite element method for the nonlinear thermistor system and
present optimal error estimates in both $L^2$ and $H^1$ norms
without any stepsize restrictions.  In this method, the standard Crank--Nicolson
scheme is applied for the linear term in the temperature equation and
an extrapolation approximation is used for the nonlinear electric conductivity.
At each time step, one only needs to solve two uncoupled linear systems.
The main idea of our aprooach is to split the error function into
two parts, the spatially discrete error and the temporally discrete error,
by introducing a corresponding time-discrete (elliptic) system.
The former arises from the Galerkin FEM
discretization for the time-discrete equations and depends only upon
the spatial mesh size $h$ (independent of the time-step size $\tau$).
If a suitable regularity of the solution to the
time-discrete equations has been proved, the numerical solution can be bounded by
\begin{equation}
\| U_h^n - R_hU^n \|_{L^\infty} \le C h^{-d/2} \| U_h^n - R_hU^n \|_{L^2} \le
C h^{-d/2} h^{r+1}
\end{equation}
without any time-step condition,
where $U^n$ is the solution of the time-discrete equations.
More important is that our approach is applicable for more general
nonlinear parabolic PDEs and many other time discretizations
to obtain unconditional convergence and optimal error estimates.


The rest of the paper is organized as follows. In Section 2, we
present the uncoupled and linearized Crank--Nicolson scheme with a linear
Galerkin finite element approximation in the spatial direction and state our
main results. After introducing the corresponding time-discrete
system, we provide in Section 3 a priori estimates and
optimal error estimates for the time-discrete solution, which imply
the suitable regularity of the time-discrete solution. With the
regularity obtained,  in Section 4 we present optimal error estimates of the
fully discrete Galerkin finite element solution in both the $L^2$ norm
and the $H^1$ norm without any time-step conditions.
Numerical results are presented in Section 5 to confirm
our theoretical analysis.

\section{The main result}
\setcounter{equation}{0}
Let $\Omega$ be a bounded, smooth and convex domain in $\R^d$
($d=2,3$). Let $\pi_h$ be a regular division of $\Omega$ into
triangles
$T_j$, $j=1,\cdots,M$ in $\R^2$ or tetrahedra in $\R^3$, and denote by
$h=\max_{1\leq j\leq
M}\{\mbox{diam}\,T_j\}$ the mesh size. For a triangle (or tetrahedra) $T_j$ at the
boundary, we define $\tilde T_j$ to be a triangle with one curved
side (or a tetrahedra with one curved face in $\R^3$) with the same
vertices as $T_j$, and set $D_j=\tilde T_j\backslash
T_j$. For an interior triangle, we set $\tilde T_j=T_j$ and
$D_j=\emptyset$. For a given triangular (or tetrahedral) division of
$\Omega$, we define the finite element spaces \cite{Tho}:
\begin{align*}
&V_{h}=\{v_h\in C(\overline\Omega):
v_h|_{T_j}\mbox{~is~linear~at~each~element~and~}v_h=0~\mbox{on}~
D_j\} ,\\
&S_{h}=\{v_h\in C(\overline\Omega): v_h|_{\widetilde T_j}
\mbox{~is~linear~at~each~element}\} .
\end{align*}
It follows that $V_h$ is a subspace of $H^1_0(\Omega)$ and $S_h$ is
a subspace of $H^1(\Omega)$. For any function $v\in S_h$, we define
$\Lambda_hv$ to be a
function satisfying $\Lambda_hv=0$ on $D_j$ and $\Lambda_hv=v$ on
$T_j$. We further define $\widetilde\Pi_h: C_0(\overline\Omega)\rightarrow
S_h$ to be the Lagrangian interpolation operator and set
$\Pi_h=\Lambda_h\widetilde\Pi_h$. Clearly, $\Pi_h$ is a projection
operator from $C_0(\overline\Omega)$ onto $V_h$.

Let $0=t_0<t_1<\cdots<t_N=T$ be a
uniform partition of the time interval $[0,T]$ with $t_n=n\tau$ and
let
\begin{align}
u^n = u(x,t_n), \quad \phi^n = \phi(x,t_n) \, \qquad\mbox{for}\quad
n=0,1,\cdots,N.
\end{align}
For any sequence of functions $\{ f^n \}_{n=0}^N$, we define
\begin{align}
& D_\tau f^{n+1}=\frac{f^{n+1}-f^n}{\tau},\qquad  \widehat
f^{n+1/2}=(3f^{n}-f^{n-1})/2 , \\
& \bar f^{n+1/2} = \frac{1}{2} ( f^n + f^{n+1}) ,
\label{dfnsjlwejriow}
\end{align}
for $n=1,2,\cdots,N-1$.

For the simplicity of notations, we denote by $C$ a generic positive
constant and by $\epsilon$ a generic small positive constant,
which depend solely upon the physical parameters of the problem and
independent of $\tau$, $h$ and $n$.
We assume that $g(\cdot,t)\in H^1(\Omega)$ is given for each fixed
$t\geq 0$.

We propose an uncoupled and linearized Crank--Nicolson Galerkin finite element
method to solve the system (\ref{e-heat-1})-(\ref{BC}), which seeks
$U^{n+1}_h\in V_h$ and $\Phi^{n+1/2}_h\in g^{n-1/2}+V_h$, $n=0,1,\cdots,N-1$, such that
\begin{align}
&\Big(\sigma(\widehat U_h^{n+ 1/2})\nabla\Phi_h^{n+ 1/2},\nabla
\varphi\Big)=0,  \quad\forall~ \varphi \in V_h,
\label{TDFEM-heat-1}
\\
&\Big(D_\tau U^{n+1}_h,v\Big)+\Big(\nabla \overline U_h^{n+ 1/2},\nabla
v\Big)=\Big(\sigma(\widehat U_h^{n+ 1/2})|\nabla\Phi_h^{n+ 1/2}|^2,v\Big),
\quad\forall~ v \in V_h
\label{TDFEM-heat-2}
\end{align}
where a standard extrapolation \cite{DFJ} is used to approximate the
nonlinear electric conductivity for $n>0$.

At the initial time steps, we choose $U^0_h=\Pi_hu_0$ and
let $\Phi^0_h$ be the Galerkin solution to the potential equation
\begin{align}
&\Big(\sigma(u_0)\nabla\Phi_h^0,\nabla \varphi\Big)=0 ,\quad
\forall~\varphi\in V_h.
\end{align}
and $\widehat U_h^{1/2}$ can be calculated either by a semi-implict Euler scheme
\begin{align}
\Big(\frac{\widehat U^{1/2}_h - u^0}{\tau/2}, v\Big)
+ \Big(\nabla \widehat U_h^{1/2}, \nabla v\Big)
=\Big(\sigma(u_0)|\nabla\Phi_h^0|^2, v\Big),
\quad \forall~ v \in V_h
\label{U-1/2-1}
\end{align}
or by an explicit Euler scheme
\begin{align}
\Big(\frac{\widehat U^{1/2}_h - u^0}{\tau/2}, v\Big) + \Big(\nabla u_0,\nabla v\Big)
=\Big(\sigma(u_0)|\nabla\Phi_h^0|^2, v\Big),
\quad \forall~ v \in V_h
\label{U-1/2-2}
\end{align}

By the classical finite element theory for elliptic equations and for interpolation,
we have
\begin{align}\label{inisteperrestimates}
\|U^0_h-u_0\|_{L^2}+\|\Phi_h^0-\phi^0\|_{L^{12/5}}\leq Ch^2  .
\end{align}

Here we assume that the solution of the initial/boundary
value problem (\ref{e-heat-1})-(\ref{BC}) exists and satisfies
\begin{align}
\label{StrongSOlEST-1}
&\|u_0\|_{H^2}+\|u\|_{L^\infty((0,T);H^2)}+\|u_t\|_{L^\infty((0,T);H^2)}  \nn\\
&+\|u_{tt}\|_{L^\infty((0,T);H^1)}
+\|u_{tt}\|_{L^2((0,T);H^2)}+
\|u_{ttt}\|_{L^2((0,T);L^2)} \le C  ,
\end{align}
and
\begin{align}
&\label{StrongSOlEST-2}
\|\phi\|_{L^\infty((0,T);W^{2,12/5})}
+\|\phi_t\|_{L^\infty((0,T);W^{1,6})}
+\|\nabla\phi\|_{L^\infty((0,T);L^\infty)} \nn \\
& + \|g\|_{L^\infty((0,T);W^{2,12/5})}
+\|\nabla g\|_{L^\infty((0,T);L^\infty)}\leq C .
\end{align}
The emphasis of this paper is on the unconditionally optimal error analysis.
The above regularity assumptions may possibly be weakened for the analysis below.
We present our main result in the following theorem.
The proof will be given in Sections 3-4.

\begin{theorem}\label{ErrestFEMSol}
{\it
Suppose that the system \refe{e-heat-1}-\refe{e-heat-2} with the
initial and boundary conditions \refe{BC} has a unique solution $(u,
\phi)$ satisfying \refe{StrongSOlEST-1}-\refe{StrongSOlEST-2}.
Then the finite element system
(\ref{TDFEM-heat-1})-(\ref{U-1/2-1}) admits a unique solution
$(U^n_h, \, \Phi^{n-1/2}_h)$, $n=1,\cdots,N$, such that
\begin{align}\label{optimalL2est}
&\max_{1\leq n\leq N}\|U^{n}_h-u(\cdot,t_{n})\|_{L^2}
+\max_{1\leq n\leq N}\|\Phi^{n-1/2}_h-\phi(\cdot,t_{n-1/2})\|_{L^{12/5}}
\leq C (\tau^2+h^2)
,\\[5pt]
&\max_{1\leq n\leq N}\|U^{n}_h-u(\cdot,t_{n})\|_{H^1}
+\max_{1\leq n\leq N}\|\Phi^{n-1/2}_h-\phi(\cdot,t_{n-1/2})\|_{W^{1,12/5}}
\leq C (\tau^2+h).
\end{align}}
\end{theorem}
\medskip

To prove Theorem \ref{ErrestFEMSol}, we introduce a time-discrete system of
equations:
\begin{align}
&\nabla\cdot\big(\sigma(\widehat U^{n+ 1/2})\nabla\Phi^{n+ 1/2}\big)=0,
\qquad\qquad\qquad\qquad~~~\mbox{for $n=0,1,\cdots$}\label{TD-heat-1}
\\
&D_\tau U^{n+1}-\Delta \overline U^{n+ 1/2}
=\sigma(\widehat U^{n+ 1/2})|\nabla\Phi^{n+1/2}|^2, \qquad~\mbox{for
$n=0,1,\cdots$}\label{TD-heat-2}
\end{align}
subject to the boundary/initial conditions
\begin{align}
\label{TDBC}
\begin{array}{ll}
U^n(x)=0,\quad \Phi^{n+1/2}(x)=g(x,t_{n+1/2})~~
&\mbox{for}~~x\in\partial\Omega,~~1\leq n\leq N,
\\[3pt]
U^0(x)=u_0(x)~ &\mbox{for}~~x\in\Omega
\\[5pt]
\displaystyle\frac{\widehat U^{1/2} - U^0}{\tau/2} - \Delta \widehat U^{1/2}
= \sigma(u_0)|\nabla\phi^0|^2 ~ &\mbox{for}~~x\in\Omega\, ,
\\[5pt]
\widehat U^{1/2}=0,~~
&\mbox{for}~~x\in\partial\Omega,
\end{array}
\end{align}
where $\phi^0$ is the solution to the elliptic equation
\begin{align}
\label{iniphi0}
\left\{
\begin{array}{ll}
\nabla\cdot(\sigma(u_0)\nabla\phi^0)=0\quad\mbox{in}
~~\Omega,\\[5pt]
\phi^0(x)=g(x,0)\quad\mbox{for}~~x\in\partial\Omega .
\end{array}\right.
\end{align}
With the solution of the time-discrete system $(U^n, \Phi^{n-1/2})$,
we have the following error splitting:
\begin{align}
&\|U^n_h- u^n\| \leq \| e^n \|
 + \| U^n- U_h^n \|
\nn \\
&\|\Phi_h^{n-1/2}- \phi^{n- 1/2}\| \leq \| \eta^{n- 1/2} \| +
\| \Phi^{n-1/2} - \Phi_h^{n- 1/2} \|
\nn
\end{align}
where
\begin{align}
& e^n = U^n - u^n,
\qquad \eta^{n-1/2} = \Phi^{n-1/2} - \phi^{n-1/2}
\, . \nn
\end{align}

Note that the fully discrete system
(\ref{TDFEM-heat-1})-(\ref{TDFEM-heat-2}) can be viewed as the
spatial discretization of the elliptic system
(\ref{TD-heat-1})-(\ref{TD-heat-2}). The
key issue is to prove the regularity of the solution to the
time-discrete equations (\ref{TD-heat-1})-(\ref{TD-heat-2}) required
in the error estimates of the Galerkin finite element method.
We present the estimates of the error functions $(e^{n}, \eta^{n- 1/2})$ and
$(U^n-U_h^n, \, \Phi^{n-1/2}-\Phi_h^{n-1/2})$ in Section 3 and Section 4, respectively.

Moreover, the error estimates given in the above theorem for
$\Phi$ are defined in the time level $t_{n-1/2}$.
To get the solution at the time level $t_{n}$, we define
$$
\Phi_h^n: = \frac{1}{2} (\Phi_h^{n+1/2}+\Phi_h^{n-1/2})
\, .
$$
By the above theorem, we see that
\begin{align}
&\max_{1\leq n\leq N}\|\Phi^n_h-\phi(\cdot,t_n)\|_{L^{12/5}}\leq
C(\tau^2+h^2)
,\\[5pt]%
& \max_{1\leq n\leq N}\|\Phi^n_h-\phi(\cdot,t_n)\|_{W^{1,12/5}}
\leq C(\tau^2+h).
\end{align}

The following lemma can be proved by noting the definition (\ref{dfnsjlwejriow})
and using a triangular inequality.

\begin{lemma}
\label{lemma1}
{\it Let $\{ v^n \}_{n=0}^N$ be a sequence of functions on $\Omega$.
Then for any norm $\| \cdot \|$,
\begin{equation}
\tau\|  v^n \| \leq 2 \tau\sum_{m=1}^n \| \overline v^{m-1/2} \|
+ \tau \| v^0 \|
\leq 2 \sqrt{T}\sqrt{ \sum_{m=1}^n \tau \| \overline v^{m-1/2} \|^2 }  + \tau \| v^0 \|
\, .
\end{equation}
}
\end{lemma}

\section{Temporal error analysis}\label{dejlagjs}
\setcounter{equation}{0}
In this section, we prove the existence and uniqueness of the solution
of the time-discrete system (\ref{TD-heat-1})-(\ref{TDBC}) and
establish error bounds for $(e^n, \eta^{n-1/2})$.

\begin{theorem}\label{ErrestDisSol}
{\it Suppose that the system \refe{e-heat-1}-\refe{BC} has a unique
solution $(u, \phi)$ satisfying \refe{StrongSOlEST-1}-\refe{StrongSOlEST-2}. Then
the time-discrete system (\ref{TD-heat-1})-(\ref{TDBC}) admits a
unique solution $(U^n, \Phi^{n-1/2})$ such that
\begin{align}
&\max_{1\leq n\leq N}\| U^n\|_{H^2}
+ \max_{1\leq n\leq N}\| D_\tau U^n\|_{H^2}+ \| D_\tau \widehat U^{1/2}\|_{H^2}
\leq C, \label{ErrestDisSol211}\\
&\max_{1\leq n\leq N}\| \Phi^{n-1/2}\|_{W^{2,12/5}} +\max_{1\leq
n\leq N}\|\nabla\Phi^{n-1/2}\|_{L^p}\leq C, \quad \forall~1\leq
p<\infty,
\label{ErrestDisSol222}
\end{align}
and
\begin{align}\label{TDErrEstLemm}
\begin{array}{ll}
&\displaystyle\max_{1\leq n\leq N}\| e^{n}\|_{H^1}
+\max_{1\leq n\leq N}\| \eta^{n-1/2}\|_{W^{1,12/5}} \leq C \tau^2.
\end{array}
\end{align}
}
\end{theorem}

\noindent{\it Proof}~~~ The existence and uniqueness of solution to the linear partial differential equations \refe{TD-heat-1}-\refe{TDBC} is obvious.
In the following, we only prove the estimates (\ref{ErrestDisSol211})-(\ref{TDErrEstLemm}).

Since $U^0=u^0$, the error functions $e^{n+1}$
and $\eta^{n+1/2}$, $0\leq n\leq N-1$, satisfy
\begin{align}
&-\nabla\cdot\big[\sigma(\widehat  U^{n+1/2})\nabla\eta^{n+ 1/2}\big]
=\nabla\cdot\big[(\sigma(\widehat u^{n+1/2})
-\sigma(\widehat  U^{n+1/2}))\nabla \phi^{n+1/2}\big]
+\nabla\cdot R_{\phi}^{n+1/2},
\label{TDerr-heat-1}
\\[8pt]
&D_\tau e^{n+1}-\Delta \overline e^{n+ 1/2}=(\sigma( \widehat U^{n+1/2}) -
\sigma(\widehat u^{n+1/2}))
|\nabla \phi^{n+1/2}|^2\label{TDerr-heat-2} \\
&~~~~~~~~~~~~~~~~~~~~~~~~~+\sigma(\widehat U^{n+1/2})(\nabla\phi^{n+1/2}
+\nabla \Phi^{n+1/2})\cdot\nabla \eta^{n+1/2}+R_u^{n+1},\nonumber
\end{align}
and
\begin{align}
\frac{2}{\tau} \widehat e^{1/2} - \Delta \widehat e^{1/2} =R_u^{1/2} ,
\label{e-1/2}
\end{align}
where
\begin{align*}
&R_u^{1/2}=
u_t |_{t=0} - \frac{u^{1/2} - u^0}{\tau/2}+\sigma(u^0)|\nabla\phi^0|^2
-\sigma(u^{1/2})|\nabla\phi^{1/2}|^2
\, \\
&R_u^{n+1} = \frac{\partial u}{\partial t}\Big|_{t=t_{n+1/2}}\!\!\!\!\!
- D_\tau u^{n+1}+\Delta (\overline
u^{n+1/2}-u^{n+1/2})+(\sigma(u^{n+1/2})
-\sigma(\widehat u^{n+1/2}))|\nabla\phi^{n+1/2}|^2
\, ,\\
&R_{\phi}^{n+1/2} = (\sigma(u^{n+1/2})-\sigma(\widehat
u^{n+1/2}))\nabla\phi^{n+1/2}
\end{align*}
are the truncation errors.
With the regularity given in \refe{StrongSOlEST-1}-\refe{StrongSOlEST-2},
we have the following estimates for the truncation errors:
\begin{align}
& \max\limits_{0\leq n\leq N-1}\|R_{\phi}^{n+1/2}\|_{L^6}
+ \left ( \sum_{n=0}^{N-1} \tau \|R_u^{n+1}\|_{L^2}^2 \right )^{\frac{1}{2}}
+ \tau  \| R^{1/2}_{u} \|_{L^6}
\leq C\tau^2 . \label{truncerr}
\end{align}

To prove \refe{ErrestDisSol211}-\refe{TDErrEstLemm},
first we study the error $\widehat e^{1/2}$.
Multiplying the equation \refe{e-1/2}
by $|\widehat e^{1/2}|^4\widehat e^{1/2}$ and integrating it over $\Omega$, we get
\begin{align*}
\| \widehat e^{1/2} \|_{L^6} \leq C \tau^2
\end{align*}
which further shows that $\tau\|\widehat e^{1/2}\|_{H^2}\leq C\tau^2$. Since $H^2\hookrightarrow C^\alpha$, it follows that $\widehat U^{1/2} \in C^{\alpha}(\overline\Omega)$.
By applying Schauder's $W^{1,p}$ estimate
\cite{ADN, ChenYZ, Simader}
to the equation \refe{TDerr-heat-1} with $n=0$,
we derive that
$$
 \| \eta^{1/2} \|_{W^{1,6}} \le C \| \widehat e^{1/2} \|_{L^6}
+ C \| R_{\phi}^{1/2} \|_{L^6} \le C \tau^2
\, .
$$
By noting the fact $e^0=0$, from the equation
\refe{TDerr-heat-2} with $n=0$, we get
\begin{align*}
& \| \nabla e^1 \|_{L^2}^2 + \tau \|\Delta e^1 \|_{L^2}^2
\le \tau (\| \widehat e^{1/2} \|_{L^2}^2 + \| (\nabla\phi^{n+1/2}
+\nabla \Phi^{n+1/2})\nabla \eta^{1/2} \|_{L^2}^2 +
\| R_u^1 \|_{L^2}^2 ) \le C \tau^4 .
\end{align*}
To conclude, we have
\begin{align}
&\|e^1\|_{H^1}+ \|\eta^{1/2}\|_{W^{1,6}}+
\tau^{1/2}\|e^1\|_{H^2}
\leq C_0\tau^2 .
\label{i-error}
\end{align}

Secondly, we present $L^2$ error estimates for the solution of
(\ref{TDerr-heat-1})-(\ref{TDerr-heat-2}). Multiplying
(\ref{TDerr-heat-1}) by $\eta^{n+1/2}$ and integrating the result
over $\Omega$, we obtain
\begin{equation}
\|\eta^{n+1/2}\|_{H^1}\leq C\|\widehat e^{n+1/2} \|_{L^2}
+ C \tau^2.
\label{TDErrEq2est}
\end{equation}
Again, multiplying (\ref{TDerr-heat-2}) by $\overline e^{n+1/2}$ and
integrating it over $\Omega$ give
\begin{align*}
&\frac{1}{2} D_\tau\biggl(\|e^{n+1}\|_{L^2}^2\biggl)
+\|\nabla \overline e^{n+1/2}\|_{L^2}^2
 \\
&\leq C\|\widehat e^{n+1/2}\|_{L^2} \|\overline e^{n+1/2} \|_{L^6} \|\nabla
\phi^{n+1/2} \|_{L^6}^2  + \|R_u^{n+1}\|_{L^2}
\|\overline e^{n+1/2}\|_{L^2}
\\
& ~~~ + \big(\sigma(\widehat U^{n+1/2})(\nabla\phi^{n+1/2} +\nabla
\Phi^{n+1/2}) \overline e^{n+1/2}, \, \nabla  \eta^{n+1/2}\big).
\end{align*}
Using (\ref{TD-heat-2}) and integrating by parts,
\begin{align*}
& |\big(\sigma(\widehat U^{n+1/2})(\nabla\phi^{n+1/2} +\nabla
\Phi^{n+1/2})\overline e^{n+1/2}, \,
\nabla \eta^{n+1/2}\big) |
\\
&  \le |\big(\sigma(\widehat U^{n+1/2})\overline e^{n+1/2}\nabla\phi^{n+1/2}, \,
\nabla \eta^{n+1/2}\big)|
\\
&  \quad + | \big( \nabla \cdot (\sigma(\widehat U^{n+1/2}) \nabla
\Phi^{n+1/2})\overline e^{n+1/2} + \sigma(\widehat U^{n+1/2}) \nabla
\Phi^{n+1/2}
\cdot \nabla \overline e^{n+1/2}, \, \eta^{n+1/2} \big)|
\\
&=|\big(\sigma(\widehat U^{n+1/2})\overline e^{n+1/2}\nabla\phi^{n+1/2}, \,
\nabla \eta^{n+1/2}\big)|
\\
& + |\big( \sigma(\widehat U^{n+1/2}) \nabla \eta^{n+1/2}
\cdot \nabla \overline e^{n+1/2}, \, \eta^{n+1/2} \big)
+\big( \sigma(\widehat U^{n+1/2}) \nabla \phi^{n+1/2}
\cdot \nabla \overline e^{n+1/2}, \, \eta^{n+1/2} \big)|
\\
& \leq C \| \overline e^{n+1/2} \|_{L^2} \|\nabla\eta^{n+1/2} \|_{L^2} + C
\|\nabla \overline e^{n+1/2} \|_{L^2}
\big(\|\nabla \eta^{n+1/2}\|_{L^2} \|\eta^{n+1/2}\|_{L^\infty}
+ \|\eta^{n+1/2}\|_{L^2}\big) \,
\end{align*}
Applying the maximum principle to the elliptic equation
(\ref{TD-heat-1}), we obtain $\| \Phi^{n+1/2} \|_{L^{\infty}} \leq
C$ and so $\| \eta^{n+1/2}  \|_{L^{\infty}} \leq C$ for
$1\leq n\leq N-1$.
It follows that
\begin{align*}
& \frac{1}{2} D_\tau\biggl(\|e^{n+1}\|_{L^2}^2\biggl)+\frac{1}{2}\|\nabla
\overline e^{n+1/2} \|_{L^2}^2 \\
&\leq C\|\widehat e^{n+1/2}\|_{L^2}^2+C\| \overline e^{n+1/2} \|_{L^2}^2
+C\| \eta^{n+1/2} \|_{H^1}^2+C\|R_u^{n+1} \|_{L^2}^2 ,
\\
& \leq C \left ( \| e^{n+1} \|_{L^2}^2 + \| e^n \|_{L^2}^2
+ \| e^{n-1} \|_{L^2}^2 \right ) + C \tau^4
\end{align*}
where we have noted \refe{TDErrEq2est} and used the inequality
$\| \cdot \|_{L^6} \le C \| \cdot \|_{H^1}$.
By applying Gronwall's inequality to the above inequality,  with
(\ref{TDErrEq2est}) we see that there exists a positive constant $\tau_1>0$ such that when $\tau<\tau_1$, we have
\begin{align}\label{fdaklgteahiegawl}
&\max_{1\leq n\leq N}\| e^{n+1}\|_{L^2}^2+\max_{1\leq n\leq N-1}\|
\eta^{n+1/2}\|_{H^1}^2+\sum_{n=0}^{N-1}\| \overline
e^{n+1/2}\|_{H^1}^2 \tau \leq
C \tau^4 .
\end{align}

Finally, we study the regularity in \refe{ErrestDisSol211}-\refe{ErrestDisSol222} and the estimate for $\|\eta^{n+1}\|_{W^{1,12/5}}$.
Note that the above estimate implies that
\begin{align}
& \|D_\tau e^{n+1}\|_{L^2}^2
+  \sum_{m=0}^{n} \tau \| D_\tau \overline e ^{m+1/2} \|_{H^1}^2 \le C \tau^2,
\label{tx-norm} \\
&\|D_\tau U^{n+1}\|_{L^2}\leq
\|D_\tau u^{n+1}\|_{L^2}+\|D_\tau e^{n+1}\|_{L^2}\leq C
\nn
\end{align}
for $0 \le n \le N-1$.
Regarding (\ref{TDerr-heat-2}) as an elliptic equation and
applying the $H^2$ estimate \cite{ChenYZ}, with (\ref{tx-norm}) we obtain
\begin{align}
\|\overline e^{n+1/2}\|_{H^2}
& \leq
C \left ( \| D_\tau e^{n+1} \|_{L^2}
+ \| \widehat e^{n+1/2} \|_{L^2}
+ \|\nabla\eta^{n+1/2}\|_{L^2}+
\|\nabla\eta^{n+1/2}\|_{L^4}^2 + \| R_u^{n+1} \|_{L^2} \right )
\nn \\
&\leq
C\tau+C\|\nabla\eta^{n+1/2}\|_{L^2}^{1/2}\|\nabla\eta^{n+1/2}\|_{L^6}^{3/2}
\nn\\
&\leq C_1\tau(1+\|\nabla\eta^{n+1/2}\|_{L^6}^{3/2}),
\label{H2estintermsPhi}
\end{align}
for $n=1,\cdots,N-1$.

Now we prove a primary estimate
\begin{equation}
\|\nabla\eta^{n+1/2}\|_{L^6}\leq 1
\label{ind}
\end{equation}
for $0\leq n\leq N-1$, by mathematical induction.
It is easy to see from
(\ref{i-error}) that \refe{ind} holds for $n=0$ if $\tau<1/C_0$.
We assume that \refe{ind} holds for $0\leq n\leq k$.
Then from (\ref{H2estintermsPhi}) and \refe{i-error} we get
$$
\|\overline e^{n+1/2}\|_{H^2} \leq (2C_1+C_0\tau)\tau,~~\mbox{for}~~ 0\leq n\leq k
$$
and by Lemma 2.1,
$$
\|e^{n+1}\|_{H^2} \le 2 \sum_{m=0}^n \| \overline e^{m+1/2} \|_{H^2}
+ \| e^0 \|_{H^2}  \leq (2C_1+C_0\tau)T ,~~\mbox{for}~~1\leq n\leq k.
$$
Hence,
$$
\|U^{n+1}\|_{H^2} \leq \|u^{n+1}\|_{H^2} + (2C_1+C_0\tau)T\leq C_3,
~~\mbox{for}~~1\leq n\leq k.$$
Since $H^2\hookrightarrow C^\alpha$ in $\R^d$ ($d=2,3$), we have
$$
\|\widehat U^{k+3/2}\|_{C^{\alpha}}
\leq C\| U^{k+1}\|_{H^2}+C\| U^k\|_{H^2} \leq C_4.
$$
With the H\"{o}lder regularity of $\sigma(\widehat U^{k+3/2})$, by applying the
$W^{1,p}$ estimate \cite{Simader} to (\ref{TDerr-heat-1}) for $n=k+1$,
we obtain
\begin{align}
\|\nabla\eta^{k+3/2}\|_{L^6}
&\leq  C_5 \| (\sigma(\widehat u^{k+3/2}) - \sigma(\widehat U^{k+3/2}))
\nabla \phi^{k+3/2} \|_{L^6} +C_5\|R_{\phi}^{k+1/2}\|_{L^6} \nn \\
&\leq  C_6\| \widehat  e^{k+3/2}\|_{L^6}  +C_6\tau^2  \nn \\
& \leq C_7\| \widehat  e^{k+3/2}\|_{H^1} +C_6\tau^2
\label{eta-6}
\\
& \leq C_8 \left ( \| e^{k+1}\|_{H^1} + \| e^k\|_{H^1} \right) +C_6\tau^2  \nn \\
& \leq C_9 \sum_{m=0}^{k+1}\| \overline e^{m+1/2}\|_{H^1} + C_6 \tau^2  \leq
C_{10}\tau,
\nn
\end{align}
where we have used Lemma 2.1 and \refe{fdaklgteahiegawl}.

By choosing $\tau<1/\max\{C_0,C_{10}\}$, we get $\|\nabla\eta^{k+3/2}\|_{L^6}
\leq 1$ and we complete the induction. Thus, we have
proved that \refe{ind} holds for $0\leq n\leq N-1$, which together with
\refe{H2estintermsPhi}  implies that $\max_{1\leq n\leq N-1}\|\overline e^{n+1/2}\|_{H^2}
\leq C\tau$. By using Lemma 2.1 again, we obtain
\begin{align}\label{H2reguTDsol}
\max_{1\leq n\leq N}\|U^n\|_{H^2}\leq C \, .
\end{align}

Since $H^2\hookrightarrow C^\alpha$ in $\R^d$ ($d=2,3$), with the
H\"{o}lder continuity of $\widehat U^{n+1/2}$, we apply
the $W^{1,p}$ estimate \cite{Simader} to (\ref{TD-heat-1}) and derive that
\begin{align}\label{estphiCalpha}
\max_{1\leq n\leq N-1}\|\nabla\Phi^{n+1/2}\|_{L^p}\leq C_p
,\quad\forall~1\leq p<\infty
\end{align}
where we have noted $\| g \|_{W^{1,p}} \le C$.
With the estimates (\ref{H2reguTDsol})-(\ref{estphiCalpha}), we can
perform the $W^{2,p}$ estimate (with $p=12/5 $)
\cite{ADN, ChenYZ} for the elliptic equation (\ref{TD-heat-1}) to
obtain
\begin{align}\label{estphiH2}
\max_{1\leq n\leq N-1}\|\Phi^{n+1/2}\|_{W^{2,12/5}}\leq C
\end{align}
where we have also noted $\| g \|_{W^{2,12/5}} \le C$.


With the above estimates, multiplying (\ref{TDerr-heat-2})
by $-\Delta\overline e^{n+1/2}$ and using the inequality
$\| \overline e^{n+1/2}\|_{H^2}\leq C\|\Delta \overline e^{n+1/2}\|_{L^2}$
(because of the boundary condition $\overline e^{n+1/2}=0$ on $\partial\Omega$),
we obtain
\begin{align*}
& D_{\tau} \| e^{n+1}\|_{H^1}^2
+ \|\overline e^{n+1/2}\|_{H^2}^2\\
&\leq C (\|\nabla\phi^{n+1/2}\|_{L^\infty}^2 +\|\nabla\Phi^{n+1/2}\|_{L^3}^2)
\|\nabla\eta^{n+1/2}\|_{L^6}^2
+ C \| \widehat e^{n+1} \|_{L^2}^2 + C \|R_u^{n+1}\|_{L^2}^2 \\
&\leq C\|\widehat e^{n+1/2}\|_{H^1}^2 + C\tau^4 ,
\end{align*}
where we have used \refe{eta-6}.
By Gronwall's inequality, we see that
\begin{align}
\max_{0\leq n\leq N-1}\| e^{n+1}\|_{H^1}^2
+\sum_{n=0}^{N-1}\tau\|\overline e^{n+1/2}\|_{H^2}^2 \leq C\tau^4
\label{e-H2}
\end{align}
and by using \refe{eta-6} again, we have
\begin{equation}
\|\nabla\eta^{n+1/2}\|_{L^{12/5}} \le C \tau^2 \, .
\end{equation}
Moreover, by \refe{e-H2} and Lemma \ref{lemma1}, we have further
\begin{align*}
\|e^{n}\|_{H^2}\leq C\tau^{-1}\biggl(\sum_{k=0}^{n-1}
\tau\|\overline e^{k+1/2}\|_{H^2}^2 \biggl)^{\frac{1}{2}}\leq C\tau ,
\end{align*}
which implies that
$\|D_\tau e^{n}\|_{H^2}\leq C\tau^{-1}\| e^{n}\|_{H^2}\leq C$.
With the regularity assumption (\ref{StrongSOlEST-1}), we see that
\begin{equation}
\| D_\tau U^n\|_{H^2} \leq\| D_\tau e^n\|_{H^2}+ \| D_\tau u^n\|_{H^2}
\leq C .
\end{equation}

So far we have proved that there exists a positive constant $\tau_0$
such that for $\tau<\tau_0$ (\ref{ErrestDisSol211})-(\ref{TDErrEstLemm}) hold.
Also we have proved
that for any $\tau>0$,
(\ref{ErrestDisSol211})-(\ref{TDErrEstLemm}) hold at the initial steps.

For $\tau\geq \tau_0$,
if we assume that
\begin{align}\label{dfsjkfowueiow}
\|\nabla \Phi^{n-1/2} \|_{L^p} + \|U^n\|_{H^2} \leq C_n
\end{align}
for $1\leq n\leq k$ (for any fixed $1 < p < \infty$), where $C_n$ is a constant dependent upon $n$.
Then, by writing the equation (\ref{TD-heat-1}) as
\begin{align}\label{fsdjkl1j324i1h5io}
\Big(1-\frac{\tau}{2}\Delta\Big)U^{k+1}
=\Big(1+\frac{\tau}{2}\Delta\Big)U^k+\sigma(\widehat
U^{k+ 1/2})|\nabla\Phi^{k+ 1/2}|^2\tau  ,
\end{align}
and applying the classical $W^{1,p}$ and $W^{2,p}$ estimates \cite{ADN, ChenYZ, Simader} of elliptic equations to
(\ref{TD-heat-2}) and (\ref{fsdjkl1j324i1h5io}), we get
$$
\|\nabla \Phi^{k+1/2} \|_{L^p} + \|U^{k+1}\|_{H^2} \leq C_{k+1}  ,
$$
where $C_{k+1}$ depends upon $C_n$, $1\leq n\leq k$. By mathematical induction, (\ref{dfsjkfowueiow}) holds for $1\leq n\leq N$. Since
$N\leq \left [ T/\tau_0 \right ]+1$, by setting
$C^* = \max_{1 \le n \le N} \{ C_n \}$ we obtain
$$
\|\nabla \Phi^{n+1/2} \|_{L^p} + \|U^{n+1}\|_{H^2} \leq C^*,
~~~n=1,2,\cdots,N \, .
$$
(\ref{ErrestDisSol211})-(\ref{TDErrEstLemm}) follow immediately.
The proof of Theorem \ref{ErrestDisSol} is complete.
\endproof\medskip

\section{Spatial error analysis}\label{dejlagjs2}
\setcounter{equation}{0}
In this section, we present error estimates
of the Galerkin finite element method for the time-discrete system
(\ref{TD-heat-1})-(\ref{TD-heat-2}).
Let $P^0_h\phi^0=g^0+\Pi_h(\phi^0-g^0)$ and
$P_h^{n-1/2}\Phi^{n-1/2}=g^{n-1/2}+\Pi_h(\Phi^{n-1/2}-g^{n-1/2})$ for
$n=1,2,\cdots,N$,
and define
$R_h:H^1_0(\Omega)\rightarrow V_h$ to be a Riesz projection
operator defined by
\begin{align*}
\!\!\!\!\!\!\!\!\!\!\!\!\!\big(\nabla (v-R_hv),\nabla
w\big)=0,\quad\mbox{for~all}~~v\in H^1_0(\Omega)~~\mbox{and}~~w\in
V_h.
\end{align*}
We summarize some basic inequalities below. The proof follows the classical finite element theory
for elliptic equations, see \cite{DFJ,Whe} for references.
\begin{align}
& \| w \|_{W^{m,p}} \le C h^{(d/p - d/q)} \| w \|_{W^{m,q}},
\quad w \in V_h,  \quad 1 \le q \le p \le \infty, \, \, m=0,1,\label{inverse}\\
&\| w \|_{W^{1,p}} \le C h^{-1} \| w \|_{L^p},\qquad 1\leq p\leq \infty,
\label{inverse2}
\\
&\|v-\Pi_hv\|_{L^p}+h\|v-\Pi_hv\|_{W^{1,p}}\leq
Ch^2|v|_{W^{2,p}},\quad p > d/2\label{ip-1}
\\[5pt]
&\|\Phi^{n-1/2}-P_h^{n-1/2}\Phi^{n-1/2}\|_{L^p}
+h\|\Phi^{n-1/2}-P_h^{n-1/2}\Phi^{n-1/2}\|_{W^{1,p}} \nn \\
& \qquad\qquad\qquad\qquad\quad~~
\leq Ch^2|\Phi^{n-1/2}-g^{n-1/2}|_{W^{2,p}}, \quad p> d/2,
\label{ip-2} \\
&\|\nabla\Pi_hv\|_{L^p} \leq
C\|v\|_{W^{1,p}},\qquad \mbox{for~all~}v\in W^{1,p}(\Omega)~~\mbox{with}~~p>d ,
\label{ip-23}
\end{align}
and
\begin{align}
&\|R_hv\|_{W^{1,p}}\leq C\|v\|_{W^{1,p}},
\quad\mbox{for all}~~v\in W^{1,p},~~1< p\leq \infty, \label{ip-333}\\
&\|v-R_hv\|_{L^p} +h\|v-R_hv\|_{W^{1,p}}\leq Ch^2 \|v\|_{W^{2,p}},
\quad\mbox{for all}~~v\in W^{2,p},~~1< p< \infty \label{ip-3} ,\\
&\|v-R_hv\|_{L^p} \leq Ch^{[(d+2p)q-dp]/(2p)}\|v\|_{W^{2,q}},
\quad dp/(d+2p)\leq q\leq p .
\label{dsfkjqqqlweio}
\end{align}

%

Let $\eta^0=\Phi_h^0-P_h^0\phi^0$ and
$$
e_h^n = U_h^n - R_h U^n, \quad \eta^{n-1/2}_h = \Phi_h^{n-1/2} -
P_h^{n-1/2} \Phi^{n-1/2} ,\quad\mbox{for}~~1\leq n\leq N.
$$
We present error estimates of the spatial discretization in the following
theorem.

\begin{theorem}\label{full-error}
{\it Suppose that the system \refe{e-heat-1}-\refe{BC} has a unique
solution $(u, \phi)$ satisfying \refe{StrongSOlEST-1}-\refe{StrongSOlEST-2}.
Then the fully-discrete finite element system
(\ref{TDFEM-heat-1})-(\ref{U-1/2-1}) admits a unique solution
$(U_h^n, \Phi^{n-1/2}_h)$, $n=1,2,\cdots,N$, such that
\begin{align}
& \max_{1\leq n\leq N}\| e^{n}_h\|_{L^2} + \max_{1\leq n\leq N}\|
\eta^{n-1/2}_h\|_{L^{12/5}}
\leq C h^2 \, ,
\label{full-error-2}\\
& \max_{1\leq n\leq N}\|\nabla e_h^{n} \|_{L^2}+\max_{1\leq n\leq
N}\|\nabla \eta_h^{n-1/2} \|_{L^{12/5}} \leq Ch \,.
\label{full-error-1}
\end{align}
}
\end{theorem}
\medskip

\noindent{\it Proof}~~~
At each time step of the scheme,
one only needs to solve two uncoupled linear discrete elliptic systems.
It is easy to see that coefficient matrices in both systems are symmetric
and positive definite.
The existence and uniqueness of the Galerkin finite element solution follows
immediately.
Since the inequality (\ref{full-error-1}) follows from (\ref{full-error-2}) via the inverse inequality (\ref{inverse2}), it suffices to prove (\ref{full-error-2}).

Let $\Phi^0=\phi^0$. The solution of the time-discrete
equations  (\ref{TD-heat-1})-(\ref{TD-heat-2}) satisfies
\begin{align}
&\big(D_\tau U^{n+1}, \, v\big) +\big( \nabla\overline U^{n+1/2}, \, \nabla v
\big) =\big(\sigma(\widehat U^{n+1/2})|\nabla\Phi^{n+1/2}|^2, \,
v\big),
\label{d-FEM-1}\\[3pt]
&\big(\sigma(\widehat U^{n+1/2}) \nabla\Phi^{n+1/2}, \, \nabla
\varphi \big)=0,\qquad\qquad\qquad\qquad\qquad~\mbox{
$n=0,1,\cdots, N$} \label{d-FEM-2}
\end{align}
for any $v, \varphi\in V_h$, and
\begin{align}
\Big( \frac{\widehat U^{1/2} - u_0}{\tau/2}, \, v\Big)
+\big( \nabla \widehat U^{1/2}, \, \nabla v \big)
=\big(\sigma(u_0)|\nabla\phi^0|^2, \, v \big), \quad v \in V_h
\label{d-FEM-1/2}
\end{align}
From the above
equations and the corresponding finite element system
(\ref{TDFEM-heat-1})-(\ref{U-1/2-1}), we find that the error
functions $e_h^{n+1}, \eta_h^{n+1/2}\in V_h$, $0\leq n\leq N$, satisfy
\begin{align}
&\big(D_\tau e^{n+1}_h,\, v\big)+\big(\nabla \overline e^{n+1/2}_h,\, \nabla v
\big)
\nn\\
&=\big(D_\tau (U^{n+1}-R_hU^{n+1}), \, v \big) +\big((\sigma(\widehat
U^{n+1/2}_h) - \sigma(\widehat U^{n+1/2}))|\nabla\Phi^{n+1/2}|^2,
\, v \big)
\nn\\
&~~~~~+2\big((\sigma(\widehat U^{n+1/2}_h)-\sigma(\widehat
U^{n+1/2}))
\nabla\Phi^{n+1/2}\cdot\nabla(\Phi^{n+1/2}_h-\Phi^{n+1/2}), \,
v\big)
\nn\\
&~~~~~+\big(\sigma(\widehat
U^{n+1/2}_h)|\nabla(\Phi^{n+1/2}_h-\Phi^{n+1/2})|^2, \, v \big)
\nn\\
&~~~~~+2\big(\sigma(\widehat U^{n+1/2})\nabla\Phi^{n+1/2}\cdot
\nabla(\Phi^{n+1/2}_h-\Phi^{n+1/2}), \,v \big)\nonumber\\
&:= \sum_{i=1}^5 I_i^{n+1/2}(v), \label{FEMErrEq1}\\[8pt]
& \big(\sigma(\widehat U^{n+1/2})\nabla\eta^{n+1/2}_h, \, \nabla \varphi \big)
= -\big((\sigma(\widehat U^{n+1/2}_h)-\sigma(\widehat
U^{n+1/2}))\nabla \Phi^{n+1/2}_h, \, \nabla \varphi \big)
\nn\\
&~~~~~~~~~~~~~~~~~~~~~~~~~~~~~~~~~~~~~~ +\big(\sigma(\widehat
U^{n+1/2})\nabla(\Phi^{n+1/2}-P_h^{n+1/2}\Phi^{n+1/2}), \, \nabla
\varphi \big), \label{FEMErrEq2223}
\end{align}
and
\begin{align}\label{sdfniowhrio}
&\Big( \frac{\widehat e^{1/2}_h - e^0_h}{\tau/2}, v \Big )
- \big ( \nabla \widehat e_h^{1/2}, \nabla v \big )\nn\\
& =
\big ( \sigma(u_0) (|\nabla\phi^0|^2 - |\nabla\Phi^0_h|^2), v \big )
+ \Big( D_{\tau/2}\widehat U^{1/2} - R_h D_{\tau/2}\widehat U^{1/2}, v \Big ),
\end{align}
for all $v, \varphi\in V_h$, where $D_{\tau/2}\widehat U^{1/2}=2(\widehat U^{1/2}- u^0)/\tau$.

First, we estimate the error functions at the initial step.
Since $\nabla\cdot(\sigma(u^0)\nabla\phi^0)=0$ in $\Omega$,
by setting $v=\widehat e^{1/2}_h$ in (\ref{sdfniowhrio}) we have
\begin{align*}
&\Big|\Big( D_{\tau/2}\widehat U^{1/2} - R_h D_{\tau/2}\widehat U^{1/2},
\widehat e^{1/2}_h \Big )\Big|
\leq C\|D_{\tau/2}\widehat U^{1/2}\|_{H^2}h^2\|\widehat e^{1/2}_h \|_{L^2}
\leq \epsilon \| \widehat e^{1/2}_h \|_{H^1} + C \epsilon^{-1} h^4
\end{align*}
and
\begin{align*}
& |\big ( \sigma(u_0) (|\nabla\phi^0|^2 - |\nabla\Phi^0_h|^2), \widehat e^{1/2}_h \big )|\\
&\le | \big(\sigma( u_0)|\nabla(\Phi^0_h-\phi^0)|^2, \, \widehat e^{1/2}_h \big)|
+ 2|\big(\sigma(u_0)\nabla\Phi^0\cdot \nabla(\Phi^0_h-\phi^0),
\,\widehat e^{1/2}_h \big)|
\\
& =  | \big(\sigma( u_0)|\nabla(\Phi^0_h-\phi^0)|^2, \, \widehat e^{1/2}_h  \big)|
+ 2|\big(\sigma(u_0)\nabla\Phi^0 (\Phi^0_h-\phi^0), \, \nabla \widehat e^{1/2}_h  \big)|
\\
& \leq C \|\widehat e^{1/2}_h\|_{L^6}
\|\nabla(\phi^0-\Phi_h^0)\|_{L^{12/5}}^2
+ C \|\Phi^0_h-\phi^0\|_{L^2} \| \nabla \widehat e^{1/2}_h \|_{L^2}
\\
&\le \epsilon \| \widehat e^{1/2}_h \|_{H^1} + C \epsilon^{-1} h^4 ,
\end{align*}
where we have used integration by parts and (\ref{inisteperrestimates}).
With the above estimates, (\ref{sdfniowhrio}) reduces to
\begin{align}
&\|\widehat e^{1/2}_h\|_{L^2}^2  \leq  \|e^0_h\|_{L^2}^2+C\tau
h^4 \leq Ch^4\, .\label{dfjksdhewbnl000}
\end{align}

Secondly, we present estimates for $\|e^{n+1}_h\|_{L^2}$ and
$\|\eta^{n+1/2}_h\|_{L^{12/5}}$ for $0\leq n\leq N-1$.
For this purpose, we take $v=\overline e^{n+1/2}_h$ in
(\ref{FEMErrEq1}) and we have
\begin{align}
I_1^{n+1/2}(\overline e^{n+1/2}_h)
& \leq
 \|\overline e^{n+1/2}_h\|_{L^2}  \|D_\tau U^{n+1}-R_h D_\tau U^{n+1}\|_{L^2}
\nn\\
&\leq
 C\|\nabla\overline e^{n+1/2}_h\|_{L^2}  \|D_\tau U^{n+1}-R_h D_\tau U^{n+1}\|_{L^2}
\nn\\
&\leq \epsilon\|\nabla \overline e^{n+1/2}_h\|_{L^2}^2
+C\epsilon^{-1} \|D_\tau U^{n+1}\|_{H^2}^2h^4,
 \nn\\
I_2^{n+1/2}(\overline e^{n+1/2}_h)
& \leq C( \|\widehat e^{n+1/2}_h \|_{L^2}
+ \| \widehat U^{n+1/2} - R_h \widehat U^{n+1/2} \|_{L^2})
\|\nabla\Phi^{n+1/2}\|_{L^6}^2\| \overline e^{n+1/2}_h\|_{L^6}
 \nn\\
& \leq C( \|\widehat e^{n+1/2}_h \|_{L^2} + h^2)
\| \nabla\overline e^{n+1/2}_h\|_{L^2}
 \nn\\
& \leq \epsilon \|\nabla\overline  e^{n+1/2}_h \|_{L^2}^2 + C\epsilon^{-1}
(\|\widehat e^{n+1/2}_h \|_{L^2}^2 + h^4),
 \nn\\
I_3^{n+1/2}(\overline e^{n+1/2}_h)
& \leq C \|\overline e^{n+1/2}_h\|_{L^6}(
\|\widehat e^{n+1/2}_h\|_{L^2}
+ \| \widehat U^{n+1/2} - R_h \widehat U^{n+1/2} \|_{L^2} )
\label{fdkhlwhwetw} \\
&~~~\, \cdot ( \|\nabla\eta^{n+1/2}_h\|_{L^6} + \|\nabla(
\Phi^{n+1/2} - P_h^{n+1/2} \Phi^{n+1/2}) \|_{L^6})
\|\nabla\Phi^{n+1/2}\|_{L^6}
\nn\\
& \leq C\|\nabla\overline e^{n+1/2}_h\|_{L^2}(\|\widehat e^{n+1/2}_h\|_{L^2}+h^2)
(h^{-d/4}\|\nabla\eta^{n+1/2}_h\|_{L^{12/5}}+ \|\Phi^{n+1/2}-g^{n+1/2}\|_{H^2})
 \nn\\
&\leq \epsilon\|\nabla\overline e^{n+1/2}_h\|_{L^2}^2
+ C\epsilon^{-1}(\|\widehat e^{n+1/2}_h\|_{L^2}^2 + h^4)
(h^{-d/2} \|\nabla \eta^{n+1/2}_h\|_{L^{12/5}}^2 +C) ,
\nn \\
I_4^{n+1/2}(\overline  e^{n+1/2}_h)
& \leq C \|\overline e^{n+1/2}_h\|_{L^6}
\|\nabla(\Phi^{n+1/2}-\Phi_h^{n+1/2})\|_{L^{12/5}}^2
 \nn\\
& \leq C \|\nabla\overline e^{n+1/2}_h\|_{L^2}
(\|\nabla\eta^{n+1/2}_h\|_{L^{12/5}}^2+ h^2 )
\nn \\
& \leq \epsilon \|\nabla\overline e^{n+1/2}_h\|_{L^2}^2 + C\epsilon^{-1}
(\|\nabla\eta^{n+1/2}_h\|_{L^{12/5}}^4+h^4) \, \nn\\
I_5^{n+1/2}(\overline e^{n+1/2}_h)
& \leq
\|\overline e^{n+1/2}_h\|_{L^6}\|\Phi^{n+1/2}_h-\Phi^{n+1/2}\|_{L^{12/5}}
\|\nabla\Phi^{n+1/2}\|_{L^{12/5}}
\nn\\
& \leq
\epsilon\|\nabla\overline e^{n+1/2}_h\|_{L^2}^2
+C\epsilon^{-1}\|\Phi^{n+1/2}_h-\Phi^{n+1/2}\|_{L^{12/5}}^2
\,. \nn
\end{align}
where we have used (\ref{ip-23}) in the estimate
of $I_3^{n+1/2}(\overline e^{n+1/2}_h)$.

By applying the $W^{1,p}$ estimate \cite{RS} to
(\ref{FEMErrEq2223}) and using \refe{ip-2}, we obtain
\begin{align} \label{dfjksdhewbnl}
&\|\nabla\eta^{n+1/2}_h\|_{L^{12/5}} \\
& \leq C \left (
\| (\sigma(\widehat U^{n+1/2}_h)-\sigma(\widehat
U^{n+1/2}))\nabla \Phi^{n+1/2}_h \|_{L^{12/5}}
+
\| \sigma(\widehat U^{n+1/2})\nabla(\Phi^{n+1/2}-P_h^{n+1/2}\Phi^{n+1/2})
\|_{L^{12/5}} \right )
\nn\\
& \leq
C \|\widehat U^{n+1/2}_h-\widehat U^{n+1/2} \|_{L^4}
(\|\nabla\eta^{n+1/2}_h\|_{L^6}
+\|\nabla P_h^{n+1/2}\Phi^{n+1/2}\|_{L^6}) + Ch  \nn\\
&\leq C(\|\widehat e^{n+1/2}_h\|_{L^4}
+\|\widehat U^{n+1/2}-R_h\widehat U^{n+1/2}\|_{L^4})(\|\nabla\eta^{n+1/2}_h\|_{L^6}
+\|\nabla P_h^{n+1/2}\Phi^{n+1/2}\|_{L^6}) + Ch \nn\\
& \leq Ch^{-d/4}  (
\|\widehat e^{n+1/2}_h\|_{L^2}+h^2 )
( h^{-d/4} \|\nabla \eta_h^{n+1/2} \|_{L^{12/5}} +C)  + Ch
\nn \\
&\leq C \left (
h^{-d/2} \|\widehat e^{n+1/2}_h\|_{L^2} \| \nabla \eta_h^{n+1/2}
\|_{L^{12/5}}
+  h^{-d/4}\|\widehat e^{n+1/2}_h\|_{L^2}
+  h^{2-d/2} \|\nabla \eta_h^{n+1/2} \|_{L^{12/5}} + h \right )  ,
\nn
\end{align}
where we have used the inverse inequality (\ref{inverse}) and (\ref{dsfkjqqqlweio}) with $q=2$ and $p=4$. To estimate $\|\eta^{n+1/2}\|_{L^{12/5}}$, we rewrite (\ref{FEMErrEq2223}) as
\begin{align*}
\big(\sigma(\widehat U^{n+1/2})
\nabla(\Phi^{n+1/2}-\Phi_h^{n+1/2}), \, \nabla \varphi
\big)+\big((\sigma(\widehat U^{n+1/2})-\sigma(\widehat
U^{n+1/2}_h)) \nabla\Phi_h^{n+1/2}, \, \nabla \varphi \big) =0
\end{align*} and apply the Nitsche technique. Define $\psi$ as
the solution of the elliptic equation
\begin{align*}
-\nabla\cdot\big(\sigma(\widehat U^{n+1/2}) \nabla\psi\big)
=|\Phi^{n+1/2}-\Phi_h^{n+1/2}|^{2/5}(\Phi^{n+1/2}
-\Phi_h^{n+1/2})
\end{align*}
with the boundary condition $\psi=0$ on $\partial\Omega$.
The solution $\psi$ to the above elliptic equation satisfies that
$$
\|\psi\|_{W^{2,12/7}}\leq
C\|\Phi^{1/2}-\Phi_h^{1/2}\|_{L^{12/5}}^{7/5}.
$$
Since
\begin{align*}
&\|\Phi^{n+1/2}-\Phi_h^{n+1/2}\|_{L^{12/5}}^{12/5}\\
&=\big(\sigma(\widehat U^{n+1/2}) \nabla(\Phi^{n+1/2}-\Phi_h^{n+1/2}),
\, \nabla \psi \big)
\\
&=\big(\sigma(\widehat U^{n+ 1/2})
\nabla(\Phi^{ n+1/2}-\Phi_h^{ n+1/2}), \, \nabla (\psi_0-\Pi_h\psi)
\big)
+\big(\sigma(\widehat U^{n+ 1/2})
\nabla(\Phi^{n+1/2}-\Phi_h^{n+ 1/2}), \, \nabla \Pi_h\psi
\big)\\
&=\big(\sigma(\widehat U^{ n+1/2})
\nabla(\Phi^{n+ 1/2}-\Phi_h^{n+ 1/2}), \, \nabla (\psi-\Pi_h\psi)
\big)
+\big((\sigma(\widehat U^{n+1/2}_h)-\sigma(\widehat U^{n+1/2}))
\nabla\Phi^{ n+1/2}_h, \, \nabla \Pi_h\psi \big)\\
&\leq C\|\nabla(\Phi^{n+ 1/2}-\Phi_h^{ n+ 1/2})\|_{L^{12/5}}
\|\nabla (\psi - \Pi_h\psi)\|_{L^{12/7}}\\
&~~~
+C\|\widehat U^{n+1/2}-\widehat U^{n+1/2}_h\|_{L^2}
\|\nabla\Phi^{n+1/2}_h\|_{L^4}\|\nabla\Pi_h\psi\|_{L^4}\\
&\leq
Ch\|\nabla(\Phi^{n+1/2}-\Phi_h^{n+1/2})\|_{L^{12/5}}
\|\psi\|_{W^{2,12/7}}\\
&~~~
+\|\widehat  U^{1/2}-\widehat U^{n+1/2}_h\|_{L^2}
\left (\|\nabla\eta_h^{n+1/2}\|_{L^{4}}
+\|\nabla P_h^{n+1/2}\Phi^{1/2}\|_{L^{4}} \right )  \|\psi\|_{W^{2,12/7}} \, ,\\
&\leq
C\|\psi\|_{W^{2,12/7}}\Big[
h\|\nabla(\Phi^{n+1/2}-\Phi_h^{n+1/2})\|_{L^{12/5}}
\\
&~~~
 +\|\widehat  U^{n+1/2}-\widehat U^{n+1/2}_h\|_{L^2}
\left (h^{-\frac{d}{6}}\|\nabla\eta_h^{n+1/2}\|_{L^{12/5}}
+\|\Phi^{n+1/2}\|_{W^{2,12/5}}+\|g^{n+1/2}\|_{W^{2,12/5}} \right )  \Big] \, ,
\end{align*}
we derive that
\begin{align}\label{fdkalghewglejek}
&\|\Phi^{n+1/2}-\Phi_h^{n+1/2}\|_{L^{12/5}}\\
&\leq
Ch\|\nabla(\Phi^{n+1/2}-\Phi_h^{n+1/2})\|_{L^{12/5}}
+C\|\widehat  U^{n+1/2}-\widehat U^{n+1/2}_h\|_{L^2}
\left (h^{-d/6}\|\nabla\eta_h^{n+1/2}\|_{L^{12/5}}
+C \right ) . \nn
\end{align}
By (\ref{dfjksdhewbnl000}) and (\ref{fdkhlwhwetw})-(\ref{fdkalghewglejek}), (\ref{FEMErrEq1}) reduces to
\begin{align}
& D_\tau\big(\|e^{n+1}_h\|_{L^2}^2\big)  + \| \nabla\overline e_h^{n+1/2}
\|_{L^2}^2
\label{e-1} \\
& \leq C\epsilon^{-1} \| D_{\tau} U^{n+1} \|_{H^1}^2 h^4
+C \left (   \| e^{n}_h\|_{L^2}^2 +\| e^{n-1}_h\|_{L^2}^2 \right )
\nn \\
&~~~+
C \left ( h^{-d/3}
\| \nabla \eta_h^{n+1/2} \|_{L^{12/5}}^2\|\widehat e^{n+1/2}_h\|_{L^2}^2 +
\| \nabla \eta_h^{n+1/2} \|_{L^{12/5}}^4 + h^4 \right ) \nn
\end{align}
for  $n=1,\cdots,N-1$.

Now we prove a primary estimate
\begin{equation}\label{fsdanhkolfasd}
\|e^{n}_h\|_{L^2}\leq h^{7/4} ,\quad 0\leq n\leq N
\end{equation}
by using mathematical induction.

By (\ref{inisteperrestimates}) and (\ref{ip-3}), this inequality
holds for $n=0$ if $h<h_1$ for some given positive constant $h_1$. If we assume that (\ref{fsdanhkolfasd}) holds for $0\leq n\leq  k$,
then from (\ref{dfjksdhewbnl000}) we know that $\|\widehat e^{n+1/2}_h\|_{L^2}\leq 2h^{7/4}+Ch^2$ for $0\leq n\leq k$ and
from the inequalities (\ref{dfjksdhewbnl})
we see that there exists a positive constant $h_2$ such that when $h<h_2$,
\begin{align}
&\|\nabla\eta^{n+1/2}_h\|_{L^{12/5}}\leq
Ch, \label{dfknhlahlwegqq}~~~\mbox{for}~~ 0 \leq n\leq k .
\end{align}
With the above inequalities, (\ref{e-1}) reduces to
\begin{align*}
D_\tau\big(\|e^{n+1}_h\|_{L^2}^2\big)  + \| \nabla e_h^{n+1/2}
\|_{L^2}^2
\leq C \left ( \| e^n_h\|_{L^2}^2
+\| e^{n-1}_h\|_{L^2}^2\right )
+ C  h^4
\end{align*}
for $0\leq n\leq k$,
which implies that there
exists a positive constant  $C_{11}$ satisfying
\begin{align}\label{sdfklsje}
\|e^{n+1}_h\|_{L^2} \leq C_{11}h^2
\end{align}
for $0\leq n\leq k$.
Therefore, $\|e^{k+1}_h\|_{L^2} <h^{7/4}$ if
$h< \min(h_1,h_2,1/C_{11}^4)$, which completes the induction.
Thus,
(\ref{dfknhlahlwegqq})-(\ref{sdfklsje}) also hold for all $0\leq n \leq N-1$, and from (\ref{fdkalghewglejek}) we derive that
$\|\Phi^{n+1/2}-\Phi_h^{n+1/2}\|_{L^{12/5}}\leq Ch^2$ and so
\begin{align}
\|\eta_h^{n+1/2}\|_{L^{12/5}}\leq \|\Phi^{n+1/2}-\Phi_h^{n+1/2}\|_{L^{12/5}}
+\|P_h^{n+1/2}\Phi^{n+1/2}-\Phi^{n+1/2}\|_{L^{12/5}}\leq Ch^2
\end{align}
for $0\leq n\leq N-1$.

So far we have proved that the estimate (\ref{full-error-2}) holds
if $h<h_0$, for some positive constant $h_0$. It remains to show
that
\begin{align}
& \max_{0\leq n\leq N}\| e^n_h\|_{L^2}
+ \max_{0\leq n\leq N-1}\| \eta^{n+1/2}_h\|_{L^{12/5}}
\leq C
\end{align}
for $h\geq h_0$. In fact, substituting $\varphi=\Phi_h^{n+1/2}$ in (\ref{TDFEM-heat-1}) we get
\begin{align*}
&\|\nabla\Phi^{n+1/2}_h\|_{L^2}\leq C,\quad n=1,2,\cdots,N.
\end{align*}
By the inverse inequalities, we have
\begin{align*}
&\|\nabla\Phi^{n+1/2}_h\|_{L^\infty} \leq C h_0^{-2},\quad n=1,2,\cdots,N
\end{align*}
and so
\begin{align*}
&\|\nabla\eta^{n+1/2}_h\|_{L^{12/5}} \leq C_{h_0},\quad n=1,2,\cdots,N  .
\end{align*}
where $C_{h_0}$ is a positive constant dependent upon $h_0$.
From \refe{TDFEM-heat-1} with $v=\overline U^{n+1/2}$
we get
\begin{align*}
D_\tau \big(\|U^{n+1}_h\|_{L^2}^2\big)  + \| \nabla \overline U_h^{n+1/2} \|_{L^2}^2
& \leq C \|\nabla\Phi^{n+1/2}_h\|_{L^\infty}^2 \| \overline  U^{n+1/2}_h\|_{L^2}\\
& \leq C \|\nabla\Phi^{n+1/2}_h\|_{L^\infty}^2 \| \nabla\overline  U^{n+1/2}_h\|_{L^2}\\
& \leq C \|\nabla\Phi^{n+1/2}_h\|_{L^\infty}^4+\frac{1}{2} \| \nabla\overline  U^{n+1/2}_h\|_{L^2}^2 ,
\end{align*}
which implies that
$$
\|U^{n+1}_h\|_{L^2}^2\leq \|U^n_h\|_{L^2}^2+C_{h_0}\tau
\leq\cdots\leq C_{h_0}T.
$$
This completes the proof of Theorem \ref{full-error}.
\endproof
\vskip0.1in

Theorem \ref{ErrestFEMSol} follows from Theorem \ref{ErrestDisSol},
Theorem \ref{full-error}, together with \refe{ip-2} and \refe{ip-3}.
\endproof

\section{Numerical results}
\setcounter{equation}{0}
In this section, we present two numerical examples to illustrate our
theoretical results. The computations are performed with the software FEnics.
\vskip0.1in

{\it Example 4.1.} We rewrite the system \refe{e-heat-1}-\refe{e-heat-2}
by
\begin{eqnarray}
&& \frac{\partial u}{\partial t} - \Delta u = \sigma(u) {|\nabla \phi|}^{2} + f_1,
\label{ex-1-1}
\\
&& - \nabla \cdot ( \sigma(u) \nabla \phi) = f_2,
\label{ex-1-2}
\end{eqnarray}
where $\Omega = (0,1)\times(0,1)$ and
$$
\sigma(u) = \frac{1}{1+u^{2}}+1 .
$$
The functions $f_{1}$, $f_{2}$,
and the Dirichlet boundary conditions are chosen corresponding to the exact
solution
$$
u(x,y,t) = \exp(x+y-t), \quad \phi(x,y,t) = 1 + \sin(x+y+t).
$$

A uniform triangular partition with $M+1$ nodes in each direction
is used in our computation. We solve the system by the
linearized Crank--Nicolson Galerkin method with linear elements
and quadratic elements, respectively.
To confirm our error estimates in the $L^2$ norm,
we choose $\tau = h$ for the linear FEM
and $\tau = h^{3/2}$ for the quadratic FEM.
We present the numerical results in Tables \ref{linear-L2}-\ref{quadratic-L2}.
We can see clearly from Tables \ref{linear-L2} that the $L^2$ errors
of the linear FEM  are proportional to $h^2$ and from
Table \ref{quadratic-L2}
that the $L^2$ errors of the quadratic FEM are proportional to $h^3$.
To see the errors in the $H^1$ norm,
we take $\tau = h^{1/2}$ for the linear FEM
and $\tau = h^{3/2}$ for the quadratic FEM, and
we present numerical results in Table \ref{linear-H1}-\ref{quadratic-H1}.
All these results are in good agreement with our theoretical analysis.

To show the unconditional stability,
we test the linearized Crank--Nicolson Galerkin method with linear elements,
$h=1/80$ and the large time steps
$\tau = h,5h,10h,20h$. We present numerical results in Table \ref{lineartest}.
The results show that the scheme is stable for large time steps,
although the numerical results with $\tau = 20h$ seem not very accurate.
\vskip0.1in

\begin{table}[h]
\begin{center}
\caption{$L^{2}$ errors of linear FEM with $h = \tau = 1/M$ (Example 4.1).}
\label{linear-L2}
\begin{tabular}{c|cccccc}
\hline
 &        & $\| U_h^n - u(\cdot ,t_n) \|_{L^2}$    &       &    \\
\hline
$t$ & $M=20$ & $M=40$ & $M=80$ & order  \\
\hline
1.0 & 7.8063e-05 & 1.9587e-05 & 4.9042e-06 & 2.00 \\
2.0 & 9.9117e-05 & 2.4975e-05 & 6.2605e-06 & 1.99 \\
3.0 & 8.7998e-05 & 2.2134e-05 & 5.5422e-06 & 1.99 \\
4.0 & 5.6591e-05 & 1.4242e-05 & 3.5666e-06 & 1.99 \\
\hline
&   &   $\| \Phi_h^n - \phi(\cdot ,t_n) \|_{L^2}$   &    &     \\
\hline
$t$ & $M=20$ & $M=40$ & $M=80$ & order  \\
\hline
1.0 & 7.2691e-05 & 1.7791e-05 & 4.3746e-06 & 2.03 \\
2.0 & 9.7524e-05 & 2.3836e-05 & 5.8610e-06 & 2.03 \\
3.0 & 1.3954e-04 & 3.4376e-05 & 8.4930e-06 & 2.02 \\
4.0 & 1.4342e-04 & 3.5393e-05 & 8.7511e-06 & 2.02 \\
\hline
\end{tabular}
\end{center}
\end{table}

\begin{table}[h]
\begin{center}
\vskip-0.2in
\caption{$L^{2}$ errors of quadratic FEM with $h = 1/M$ and
$\tau = h^{3/2}$ (Example 4.1).}
  \label{quadratic-L2}
\begin{tabular}{c|ccccc}
\hline
    &   & $\| U_h^n - u(\cdot ,t_n) \|_{L^2} $   &    &  \\
\hline
$t$ & $M=10$ & $M=20$ & $M=40$ & order  \\
\hline
1.0 & 8.8214e-06 & 1.2113e-06 & 1.5496e-07 & 2.92 \\
2.0 & 2.0196e-05 & 2.5321e-06 & 3.1437e-07 & 3.00 \\
3.0 & 1.7212e-05 & 2.1149e-06 & 2.5866e-07 & 3.03 \\
4.0 & 4.7866e-06 & 5.5294e-07 & 6.3392e-08 & 3.12 \\
\hline
&   &  $\| \Phi_h^n - \phi(\cdot ,t_n) \|_{L^2}$         &    &     \\
\hline
$t$ & $M=10$ & $M=20$ & $M=40$ & order  \\
\hline
1.0 & 3.4542e-05 & 4.0535e-06 & 5.0059e-07 & 3.05 \\
2.0 & 3.9065e-05 & 4.7335e-06 & 4.3953e-07 & 3.24 \\
3.0 & 4.1164e-05 & 3.2571e-06 & 6.2613e-07 & 3.02 \\
4.0 & 3.0560e-06 & 2.5956e-05 & 3.7303e-07 & 3.06 \\
\hline
\end{tabular}
\end{center}
\end{table}

\begin{table}[p]
\begin{center}
\vskip-0.2in
 \caption{$H^{1}$ errors of linear FEM with
$h = 1/M$ and $\tau = h^{1/2}$ (Example 4.1).}
  \label{linear-H1}
\begin{tabular}{c|cccc}
\hline
 &  &  $\| U_h^n - u(\cdot ,t_n) \|_{H^1}$ &      & \\
\hline
$t$ & $M=40$ & $M=80$ & $M=160$ & order \\
\hline
1.0 & 5.6024e-03 & 2.3706e-03 & 9.9037e-04 & 1.25 \\
2.0 & 3.6159e-03 & 1.8179e-03 & 7.7195e-04 & 1.11 \\
3.0 & 2.8675e-03 & 1.3254e-03 & 6.4766e-04 & 1.07 \\
4.0 & 1.7923e-03 & 8.2924e-04 & 3.5460e-04 & 1.17 \\
\hline
&    &  $\| \Phi_h^n - \phi(\cdot ,t_n) \|_{H^1}$     &   &    \\
\hline
$t$ & $M=40$ & $M=80$ & $M=160$ & order \\
\hline
1.0 & 4.1118e-03 & 2.2510e-03 & 1.1312e-03 & 0.93 \\
2.0 & 4.8396e-03 & 2.0593e-03 & 1.3143e-03 & 0.94 \\
3.0 & 5.0531e-03 & 2.7701e-03 & 1.0180e-03 & 1.16 \\
4.0 & 5.1689e-03 & 1.9872e-03 & 1.4263e-03 & 0.93 \\
\hline
\end{tabular}
\end{center}

\begin{center}
\caption{$H^{1}$ errors of quadratic FEM with $h = \tau = 1/M$ (Example 4.1).}
 \label{quadratic-H1}
\begin{tabular}{c|cccc}
\hline
  &      & $\| U_h^n - u(\cdot ,t_n) \|_{H^1}$  &     & \\
\hline
$t$ & $M=10$ & $M=20$ & $M=40$  & order  \\
\hline
1.0 & 1.6700e-03 & 3.1162e-04 & 5.9398e-05 & 2.41 \\
2.0 & 1.3279e-03 & 2.8247e-04 & 6.4458e-05 & 2.18 \\
3.0 & 1.0156e-03 & 2.2534e-04 & 5.2879e-05 & 2.13 \\
4.0 & 5.0056e-04 & 9.5814e-05 & 1.9466e-05 & 2.34 \\
\hline
&   &$\| \Phi_h^n - \phi(\cdot ,t_n) \|_{H^1}$    &   &    \\
\hline
$t$ & $M=10$ & $M=20$ & $M=40$  & order  \\
\hline
1.0 & 1.8020e-03 & 4.4199e-04 & 1.0614e-04 & 2.04 \\
2.0 & 1.6539e-03 & 4.0090e-04 & 9.4222e-05 & 2.07 \\
3.0 & 1.6141e-03 & 3.8446e-04 & 8.8813e-05 & 2.09 \\
4.0 & 1.5725e-03 & 3.7110e-04 & 8.5031e-05 & 2.10 \\
\hline
\end{tabular}
\end{center}

\begin{center}
\caption{$L^{2}$ Errors of linear FEM with $h = 1/M$ and $\tau = kh$
(Example 4.1).}
  \label{lineartest}
\begin{tabular}{c|cccc}
\hline
  &   & $\| U_h^n - u(\cdot ,t_n) \|_{L^2}$     &       &    \\
\hline
$t$ & $k=1$ & $k=5$ & $k=10$ & $k=20$ \\
\hline
1.0 & 4.9042e-06 & 3.6772e-05  & 2.2480e-04 &  1.6693e-03 \\
2.0 & 6.2605e-06 & 8.0524e-05  & 3.1973e-04 &  1.5965e-03 \\
3.0 & 5.5422e-06 & 6.6692e-05  & 2.6145e-04 &  1.1952e-03 \\
4.0 & 3.5666e-06 & 1.7407e-05  & 6.6093e-05 &  3.9786e-04 \\
\hline
&    &  $\| \Phi_h^n - \phi(\cdot ,t_n) \|_{L^2}$  &      &    \\
\hline
$t$ & $k=1$ & $k=5$ & $k=10$ & $k=20$ \\
\hline
1.0 & 4.3746e-06 & 1.4126e-04 & 5.5525e-04 & 1.3714e-03 \\
2.0 & 5.8610e-06 & 1.2561e-04 & 4.9591e-04 & 1.0407e-03 \\
3.0 & 8.4930e-06 & 1.1959e-04 & 4.6941e-04 & 8.5232e-04 \\
4.0 & 8.7511e-06 & 1.1356e-04 & 4.4541e-04 & 7.2633e-04 \\
\hline
\end{tabular}
\end{center}
\end{table}

\newpage

{\it Example 4.2.}
In the second example, we consider the
system \refe{ex-1-1}-\refe{ex-1-2} in the three-dimensional space with
the exact solution
\begin{align}
& u(x,y,z,t) = \exp(2x+y-z)(2t+\sin(t)),
\label{ex-2-1} \\
& \phi(x,y,z,t) = \sin(x-2y)\cos(z)\exp(t) \, .
\label{ex-4-2-2}
\end{align}

We use a uniform tetrahedral partition with $M+1$ nodes
in each direction (see Figure 1). The total number
of tetrahedra is $6M^{3}$ and the total number of vertices is
$(M + 1)^{3}$. We solve the system by the proposed
Crank--Nicolson Galerkin method with linear elements.
Table \ref{linear-ex2} contains the $L^2$ errors of the numerical solution
with $\tau = h$ and $h=1/10, 1/20, 1/40$.
Similarly, we can see
that the $L^{2}$ errors for both $u$ and $\phi$ are proportional to
$h^{2}$.

Previous analysis for the three-dimensional problem
often requires a stronger time-step condition
than that for the two-dimensional problem.
Finally, we test the linear Galerkin method with $h=1/80$ and large time steps
$\tau = h, 5h, 10h$. The results are presented
in Table \ref{lineartest-ex2}. Numerical results show that the scheme is
unconditionally stable.
\vskip0.1in

\begin{figure}[h]
\begin{center}
\epsfig{file=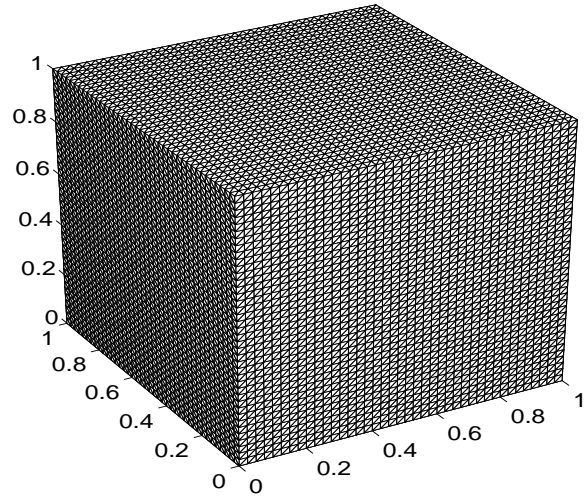,height=3 in,
width=5.0in}
\caption{The three-dimensional mesh (Example 4.2).} \label{mesh}
\end{center}
\end{figure}

\begin{table}[h]
\begin{center}
\caption{$L^{2}$ Errors of linear FEM
with  $h = \tau = 1/M$ (Example 4.2).}
  \label{linear-ex2}
\begin{tabular}{c|cccc}
\hline
 &    &  $\| U_h^n - u(\cdot ,t_n) \|_{L^2}$  &    & \\
\hline
$t$ & $M=10$ & $M=20$ & $M=40$ & order  \\
\hline
1.0 &  1.1089e-03  & 2.8319e-04  & 7.1320e-05  & 1.9793   \\
2.0 &  8.6316e-04  & 2.2523e-04  & 5.6987e-05  & 1.9605   \\
3.0 &  4.0520e-04  & 1.0626e-04  & 2.6895e-05  & 1.9566   \\
4.0 &  3.6125e-04  & 9.4243e-05  & 2.3822e-05  & 1.9613   \\
\hline
\hline
 &   & $\| \Phi_h^n - \phi(\cdot ,t_n) \|_{L^2}$    &    &     \\
\hline
$t$ & $M=10$ & $M=20$ & $M=40$ & order  \\
\hline
1.0 &  4.0129e-04  & 1.1562e-04  & 2.9779e-05  & 1.8761  \\
2.0 &  7.8577e-04  & 2.1723e-04  & 5.5611e-05  & 1.9103  \\
3.0 &  7.7231e-04  & 2.1482e-04  & 5.5087e-05  & 1.9047  \\
4.0 &  5.0533e-04  & 1.4801e-04  & 3.8404e-05  & 1.8589  \\
\hline
\end{tabular}
\end{center}
\end{table}

\begin{table}[h]
\begin{center}
\caption{$L^{2}$ Errors of linear FEM with  $h = 1/M$ and $\tau = kh$
(Example 4.2)}
  \label{lineartest-ex2}
\begin{tabular}{c|ccc}
\hline
 &    & $\| U_h^n - u(\cdot ,t_n) \|_{L^2}$   &    \\
\hline
$t$ & $k=1$ & $k=5$ & $k=10$ \\
\hline
1.0 & 7.1320e-05  & 2.5286e-04 & 2.3962e-03 \\
2.0 & 5.6987e-05  & 2.4133e-04 & 1.6463e-03 \\
3.0 & 2.6895e-05  & 1.0359e-04 & 8.5356e-04 \\
4.0 & 2.3822e-05  & 2.6147e-04 & 1.1001e-03 \\
\hline
 &    & $\| \Phi_h^n - \phi(\cdot ,t_n) \|_{L^2}$    &  \\
\hline
$t$ & $k=1$ & $k=5$ & $k=10$ \\
\hline
1.0 & 2.9779e-05 & 2.7311e-04 &  4.6730e-04 \\
2.0 & 5.5611e-05 & 1.2094e-04 &  1.3614e-03 \\
3.0 & 5.5087e-05 & 1.0514e-04 &  1.7040e-03 \\
4.0 & 3.8404e-05 & 1.0333e-04 &  1.6532e-03 \\
\hline
\end{tabular}
\end{center}
\end{table}

\newpage
\section{Conclusions}
\setcounter{equation}{0}
We have presented an uncoupled and linearized Crank--Nicolson Galerkin
finite element method for the nonlinear time-dependent thermistor equations
in the $d$-dimensional space ($d=2,3$) and
provided unconditionally optimal error estimates in both
$L^2$ and $H^1$ norms,
while existing analysis requires certain time-step restrictions.
Our numerical results confirm our analysis and show that the proposed scheme
is efficient.  Our approach presented in this paper can be extended to many
other nonlinear parabolic systems,  high-order
finite element approximations and other time discretization schemes,
while the analysis only focuses on the electric heating model
with a linear finite element method to illustrate our idea.

\end{document}